\documentclass[12pt]{article}

\usepackage{lineno,hyperref}
\usepackage{color}
\usepackage{amsfonts}
\usepackage{amsmath}
\usepackage{amssymb}
\usepackage{amsthm}
\usepackage{subeqnarray}
\usepackage{lscape,graphicx,url}
\usepackage[font=normalsize]{caption}
\usepackage{anysize}
\usepackage{proof}
\usepackage{fancyhdr}
\usepackage{latexsym}
\usepackage{array}
\usepackage{multirow}
\usepackage{textcomp}
\usepackage{optprog}
\usepackage{mathtools}
\usepackage{natbib}

\usepackage{bm}
\usepackage{array}

\usepackage{proof}
\usepackage{paralist}
\usepackage[dvipsnames]{xcolor}
\usepackage{placeins}
\usepackage{enumitem}
\usepackage{csquotes}
\usepackage{comment}
\usepackage[none]{hyphenat}

\usepackage{tikz}
\usepackage{float}

\newcolumntype{H}{>{\setbox0=\hbox\bgroup}c<{\egroup}@{}}

\begin{document}
\setlength{\bibsep}{0.0pt}

\begin{center}
{\LARGE Constructing option portfolios} \\
\vspace{0.3cm}
 J.~E.~Beasley \\
 \vspace{0.3cm}
 Mathematics, Brunel University, UK~~~john.beasley@brunel.ac.uk 
 \\
November~2026 \\
\end{center}

\section*{Abstract}
In this paper we consider the problem of constructing a portfolio composed of a set of financial options.
We present a mixed-integer linear program.
We apply our approach to index options for the S\&P~500. 

 \vspace{0.3cm}
{\bf Keywords:}   finance, index options, options, portfolio optimisation

\section{Introduction}

In this paper we consider the problem of constructing a portfolio composed of a set of financial options. We present an approach to constructing a portfolio of options,
all with the same expiry, that achieves a positive payoff at expiry with certainty (so a positive payoff irrespective of the value of the underlying at expiry).

The structure of this paper is as follows. In Section~\ref{sec2} we provide a brief introduction to options and index options, focusing on index options for the S\&P~500.
We go on to survey the literature with regard to previous work dealing with mathematical programming as applied to
 options. We state what we believe to be the contribution of this paper to the literature.

In Section~\ref{sec3} we give our mixed-integer linear program that finds (if possible) a portfolio  of options, all with the same expiry, that will collectively at expiry provide a positive payoff irrespective of the value of the underlying at expiry. 
We indicate how we can choose the values explicitly considered for the underlying which enable us to implicitly consider all possible underlying values. We go on to show how our formulation can be adapted to automatically decide the relevant investment in long and short positions given a total joint budget.
We discuss how our approach to option portfolio construction might be used in practice. 

In Section~\ref{sec4} we give our computational results when our approach is applied to 925 different problems based on S\&P~500 market data. Finally in Section~\ref{sec5} we give our conclusions.

\section{Options and literature review}
\label{sec2}

\subsection{Options and index options}

In this section we briefly consider options, focusing on the basic terminology needed to 
understand options and  index options sufficiently for the purposes of this paper. Readers less familiar with options, or wishing to learn more, are referred to textbooks 
such as~\cite{elton14} or~\cite{hull21}.

A call option on an underlying asset gives the holder of the option the right (but not the obligation) to buy the asset at
a fixed price, called the strike (or exercise) price. If this right can only be exercised at a specific (expiry) time in the future then it is known as an European call option. If it can be exercised at any time up to expiry then it is known as an American call option.

A put option is the right to sell the asset, otherwise it has the same characteristics
as a call option, so a strike price, expiry time, European or American.

\sloppy Options can be priced (valued) using an explicit mathematical formula due to \cite{black73}, sometimes referred to as  Black-Scholes-Merton  based on~\cite{merton73}. As it is an explicit mathematical formula the partial differential of this formula with respect to each and every  parameter in it can be calculated. 
These partial differentials give insight into how sensitive the price is to each parameter.
Collectively  these partial differentials are known as the \textbf{\emph{Greeks}}.  The differentials for a single option can be combined together for a portfolio of options. Achieving a portfolio with  Greek-neutrality (so portfolio differentials equal to zero) is a theme encountered in the literature.

Clearly options can apply if we have a tangible asset, e.g.~a call option to buy one tonne of a particular metal. 
Similarly
in financial markets we can have options on specific equity assets, e.g. a call option to buy 100 shares of a particular stock. 

However options also apply if we have intangible assets, such as the S\&P~500 index. As we cannot \enquote{buy} any of the index an option such as this is cash settled. So the strike price relates to the level of the index. Considering a call option, if at expiry the index is below the strike price the option is worthless. However if at expiry the index is above the strike price then the option is profitable, since the call option gives us the right to buy at a value below its market value, and so the option is cash settled. The  option holder receives a cash amount linearly related to the difference between the current index level and the strike price.
S\&P~500 index options are European options and so can only be exercised at expiry.

Options can be traded (bought/sold) over their lifetime, but typically 
only have a short life before expiry. For example the two most liquid S\&P~500 index options are: Monthly, expiring on the 
third Friday of every month and EOM, expiring on the last business day of each month. 
Exercise prices (the level of the S\&P~500 index at expiry)
are only integer multiples of 5. 


\subsection{Literature review}

Our approach to constructing a portfolio of options
is based on mathematical programming, specifically making use of a  mixed-integer linear program.
Hence in this review we focus on option papers that also make use of mathematical programming techniques such as linear/mixed-integer/integer programming.

\cite{mehta82} used integer programming to investigate buying call and put options so as to maximise expected profit.  They did not consider 
selling options, only their purchase.
They made use of Black-Scholes  pricing
both to define the  profit from an option and to calculate 
the probabilities needed to compute expected profit. Computational results were presented for five small problems.

\cite{rendleman95} used a linear program (LP) to construct an option portfolio with 
Greek-neutrality. The objective function adopted was to maximise the profit 
from the difference between the Black-Scholes price and the market price minus the transaction cost.

\cite{pap04} built on the earlier work of~\cite{rendleman95} and presented a linear programming formulation for the problem of deciding an appropriate stock (equity)
option hedging strategy, having regard to the  Greeks.
In their LP they introduced constraints to ensure Greek-neutrality. Computational results were presented for a small number of options relating  to a single company 
traded on the Stockholm Stock Exchange. 

\cite{pap04} noted that the variables in their model should be regarded as integer, rather than continuous, variables.  However in their paper achieving exact Greek neutrality 
with integer variables would be numerically difficult/impossible since their equations for doing so involve constraints forcing a weighted sum of integer variables to 
equal zero, but where the weights are given to a significant number of decimal places (three to five decimal places in their example).
\cite{horasanli08} extended the work of~\cite{pap04} and considered options on multiple assets 
and too  included any underlying assets.

\cite{gao09} presented a mixed-integer program to decide an option portfolio with risk bounds on the 
 Greeks. 
They noted that setting their risk bounds to zero meant their model led to the model of~\cite{pap04}.
 They applied their approach to the example given in~\cite{pap04}, but to do this they relax the integrality requirement.

\cite{schyns10} presented a two-period mixed-integer program, making use of scenario trees to maximise expected portfolio return at the end of the second period. Binary variables were necessary to model their Value-at-Risk constraint, making use of portfolio values at the end of each second level subtree. In order to model the value of the underlying asset over the periods considered they used
 an implied distribution approach~(\cite{schyns05}). Computational results were given for one problem involving the S\&P~500
over a one year period with rebalancing occurring only once, after   two months.

\cite{yin13} presented an approach based upon multi-stage stochastic programming 
for portfolio optimisation involving stock indices and options. Their approach involves a linear program.
They introduced risk bounds on the Greeks.
Their approach makes use of scenario trees and computational results were presented for a problem involving four market indices,

\cite{lin16} considered the problem of constructing option portfolios with Greek-neutrality and discussed the models given 
by~\cite{rendleman95, pap04, horasanli08}. They presented a linear program for producing a Greek-neutral  option portfolio with multiple underlying assets. They also presented a fuzzy goal programming model to deal with minimum transaction lot sizes. Computational results were presented for one example problem.

\cite{constantinides20} presented a linear program to find a set of option portfolios that stochastically dominate (in the second degree)  an \enquote{index trading} portfolio composed of a risk-free bond and a fund tracking a market index. Their linear program maximises expected payoff and their approach makes use of forecasts of future index values. Computational results are presented for options on the S\&P~500 index, as well as for the French CAC and German DAX indices. 

\cite{yin20} presented an approach based upon multi-stage stochastic programming 
for portfolio optimisation involving a number of different assets, but including options. Their approach makes use of scenario trees with simulated data for future option prices. Their approach entails solution of a linear program. Computational results were presented for a problem involving four market indices:  S\&P~500, FTSE~100, DAX~30 and Nikkei~225
together with  options on these indices.

\cite{post21} used a linear program to identify option portfolios that exhibit second-order stochastic  dominance over the index and hence constitute an arbitrage opportunity. They made use of probabilities for the value of the underlying at option expiry derived from a conditional multinomial distribution.
Computational results were given for S\&P~500 index options.

In a recent paper~\cite{beare25} considered the problem of constructing a portfolio of options which when combined with unit investment in the underlying asset provides a payoff that stochastically dominates investment in the underlying asset at a fixed time in the future. 

To achieve 
second-order stochastic dominance they solved a LP based on~\cite{luedtke08} and demonstrate that their LP is computationally superior to the LP given in~\cite{post21}.
To achieve first-order stochastic dominance they used a mixed-integer program involving binary variables. They presented computational results for  monthly S\&P~500 index options but found no evidence that  arbitrage opportunities  are systematically present. 

 Their approach requires use of the probabilities associated with the level of the S\&P~500 at a specified time in the future. To do this they
assumed that the S\&P~500 return distribution can be represented by a skewed generalised 
t-distribution~(\cite{theodossiou98}) with distribution parameters being calculated from current market data.

\subsection{Contribution}

The work presented in this paper is fundamentally different from literature considered above. In this paper we are concerned not with any stochastic measures, but rather with certainty. Our formulation makes no use of Greeks, probabilities (either point values or distributions) or scenarios.

\textbf{\emph{We present a mixed-integer linear program that finds (if possible) a set of options, all with the same expiry, that will collectively at expiry provide a positive payoff irrespective of  the value of the underlying at expiry.}}

To the best of our knowledge finding an  option portfolio that achieves a positive payoff with certainty in the manner described below has never been given in the literature previously.

\section{Formulation}
\label{sec3}
In this section we first present our mixed-integer linear program for constructing an option portfolio. We then indicate how we can choose the values explicitly considered for the underlying which enable us to implicitly consider all possible underlying values. We go on to indicate how the formulation can be adapted to automatically decide the relevant investment in long and short positions given a total joint budget.
Finally we discuss how our approach  can  be used in practice.

\textbf{Note here that we  assume that all of the options considered have the same expiry time}. Let:
\begin{compactitem}
\item $C$ be the set of options, where for ease of presentation we assume that any option $i \in C$ 
 can be 
 either bought or sold, so here a decision must be made as to whether to hold a long or short position in each option
 \item $E_i$ be the exercise (strike) price of option $i \in C$ 
 \item $\tau(i)$ be the type  of option $i \in C$, where we refer to a call option as an option of type 1; a put option as an option of type 2
\item $\Gamma_B$ be the initial cash budget available for long positions
\item $\Gamma_S$ be the limit on the total short position
\item $c_{ik}$ be the cost of buying one unit of option 
$i \in  C$ of type $k{=}\tau(i)$, with regard to options this cost is known as the ask price  
\item $r_{ik}$ be the revenue from selling one unit of option 
$i \in C$ of type $k{=}\tau(i)$, with regard to options this revenue is known as the bid price   
\item $\Delta =[U_j~|~j=1,\ldots,N]$ be the set of prices (values)
considered for the underlying at the chosen expiry time, where $U_j$ is the $j^{th}$ price and these prices are 
distinct and are in ascending order, so $U_1 < U_2 < \ldots < U_N$  \\
\end{compactitem}

\noindent The payoff, at expiry, associated with buying or selling any particular option $i \in C$ is as given in Table~\ref{table1}. 
In this table the notation adopted is that the first subscript for $P$ refers to the option $i$;
 the second subscript refers to the option type (so $\it{1}$ or $\it{2}$); the third subscript $j$ refers to the price of the underlying $U_j$ at expiry, and the fourth  subscript refers to whether  the option is bought ($b$) or sold ($s$).

For example the payoff associated with buying a call option is firstly the cost, so $-c_{i1}$ and then any revenue if at expiry we have $U_j \geq E_i$, since in such a case we have the right to buy an asset worth $U_j$ in the market at a price of $E_i$. This revenue is hence max$[0,U_j - E_i]$, giving the total payoff  of 
$P_{i1jb} = -c_{i1} + max[0,U_j - E_i]$ as
seen in Table~\ref{table1}.

\begin{table}[!ht]
\centering
{
\renewcommand{\tabcolsep}{1mm} 
\renewcommand{\arraystretch}{1.4} 
\begin{tabular}{|c|c|c|}
\hline
\multicolumn{1}{|l}{Option type}& 
\multicolumn{1}{|c}{Bought}& 
\multicolumn{1}{|c|}{Sold}
\\
\hline
Call, type 1, $\tau(i)=1$  & $P_{i1jb}= -c_{i1} + max[0,U_j - E_i]$ & $P_{i1js} = r_{i1} -max[0,U_j - E_i]$ \\
\hline
Put, type 2, $\tau(i)=2$  & $ P_{i2jb} =-c_{i2}  +max[0,E_i - U_j ]$ & $P_{i2js} = r_{i2} -max[0,E_i - U_j ]$\\
\hline
\end{tabular}
}
\caption{Payoffs}
\label{table1}
\end{table}

\noindent Then the decision variables are:
\begin{compactitem}
\item $x_i \geq 0$  the integer number of  units of option $i \in C$ that are bought, so held long
\item $y_i \geq 0$ the integer number of units of option $i \in  C$ 
 that are sold, so held short
\item $z_i=1$ if option $i \in C$ is bought, =0 otherwise (if option $i \in C$ is sold); here we assume that 
we need for each option to decide whether to buy or sell the option
\item $\alpha_j $ the total payoff (profit) at expiry when the value of the underlying is $U_j$~$j=1,\ldots,N$
\item $\mathcal{V}$ the minimum payoff over all values  $\alpha_j~j=1,\ldots,N$ \\
\end{compactitem}

\noindent Then the problem of determining an option portfolio is:

\begin{equation}
\mbox{maximise} ~\mathcal{V}
\label{eqobvj}
\end{equation}

\noindent subject to

\begin{equation}
\sum_{i \in C~k{=}\tau(i)} r_{ik} y_i \leq \Gamma_S
\label{eq1}
\end{equation}

\begin{equation}
\sum_{i \in C~k{=}\tau(i)} c_{ik} x_i   \leq  \Gamma_B+ \sum_{i \in C~k{=}\tau(i)} r_{ik} y_i 
\label{eq2}
\end{equation}

\begin{equation}
 x_i   \leq [( \Gamma_B + \Gamma_S)/c_{ik}]  z_i~~~~\forall i \in C~k{=}\tau(i)
\label{eq4}
\end{equation}

\begin{equation}
 y_i   \leq  (\Gamma_S/r_{ik}) (1-z_i)~~~~\forall i \in  C~k{=}\tau(i)
\label{eq3}
\end{equation}

\begin{equation}
\alpha_j = \sum_{k=1}^2 \left[ ~ \sum_{i \in C~\tau(i)=k}  (P_{ikjb} x_i 
+ P_{ikjs} y_i) ~\right] ~~~~j=1,\ldots,N
\label{eq5}
\end{equation}

\begin{equation}
\mathcal{V} \leq \alpha_j ~~~~j=1,\ldots,N
\label{eqvj}
\end{equation}

\begin{equation}
\mathcal{V}  \geq 0
\label{need1}
\end{equation}

\begin{equation}
\alpha_j \geq 0~~~~j=1,\ldots,N
\label{need2}
\end{equation}

\begin{equation}
z_i \in \{0,1\}~~~~\forall i \in C
\label{need3}
\end{equation}

\begin{equation}
x_i, y_i \geq 0, integer~~~~\forall i \in C
\label{need4}
\end{equation}

Equation~(\ref{eq1}) limits the total short position. 
Equation~(\ref{eq2}) limits the total long position, where the right-hand side is the initial budget plus any income from short sales. 

If $z_i=1$, so option $i \in C$ is bought, Equation~(\ref{eq4}) 
limits the total number of units bought to the total budget (long plus short) divided by  the cost of 
buying one unit of option $i$.
If $z_i=0$, so option $i \in C$ is not bought, Equation~(\ref{eq4}) ensures that $x_i =0$.

 If $z_i=0$, so option $i \in C$ is sold, Equation~(\ref{eq3})
limits the total number of units sold to the budget divided by  the revenue from selling one unit of option $i$.
If $z_i=1$, so option $i \in C$ is bought, Equation~(\ref{eq3}) ensures that $y_i=0$, so the option cannot be sold.

Collectively Equations~(\ref{eq4}),(\ref{eq3}) together with the binary zero-one 
nature of $z_i$ ensure that we can never have a solution in which both $x_i > 0$ and $y_i > 0$.
 If no trading is done with regard to option $i \in C$, so $x_i=y_i=0$, then  the value assigned to  $z_i$ is irrelevant.

Note that above we have 
assumed that each option $i \in C$ can be either bought or sold. If for a particular option  $i$ it cannot be bought 
we simply add  the explicit constraint that $x_i=0$. Similarly if option $i$ cannot be sold we add the explicit constraint that $y_i=0$. 

Equation~(\ref{eq5}) defines the payoff $\alpha_j$ associated with the underlying at expiry being $U_j$.
In this equation the summation is over the four payoff
terms given in Table~\ref{table1}, but where only two of these terms will be relevant, the two associated with $x_{i\tau(i)}$
 and $y_{i\tau(i)}$.

Equation~(\ref{eqvj}) ensures that $\mathcal{V}$ is less than or equal to 
each $\alpha_j$ value. Since we are maximising this means that, at optimality, the value associated with $\mathcal{V}$ must be equal to the minimum $\alpha_j$ value.

Note here that the solution $\mathcal{V}=0$ associated with  $x_i=y_i=z_i=0~\forall i \in C$ (implying $\alpha_j=0~j=1,\ldots,N$) is a feasible solution to the above formulation. Knowing that this feasible solution exists justifies the non-negativity equations, 
Equations~(\ref{need1}),(\ref{need2}). 

Equation~(\ref{need3}) defines
$z_i$ as a zero-one variable and Equation~(\ref{need4}) defines $x_i$ and $y_i$ as non-negative integer variables.

As we are maximising  any solution found will either have $\mathcal{V}=0$  or will have $\mathcal{V}>0$. 
 If $\mathcal{V}>0$ then  this implies from Equation~(\ref{eqvj}) that the 
 payoff $\alpha_j$ at each and every one of the considered underlying values $[U_1,\ldots,U_N]$ are positive, so $\alpha_j > 0~j=1,\ldots,N$.  
\textbf{\emph{In other words we will have found a set of options which achieve a positive payoff at every one of the underlying values explicitly considered.}}

With regard to a minor point it is clear that the variables $[\alpha_j]$ seen above can be eliminated by algebraic substitution. 
However we have retained them in the formulation given above for clarity of exposition.

It is trivial to introduce  into the formulation given above limits on the number of options that can be purchased or sold at the current prices.
Such limits are referred to in the options literature as market~\enquote{size}or market~\enquote{liquidity or depth}. They take the form of
 simple inequalities. For option $i \in C$ let $\delta_i$ be the upper limit on the number that can be bought and let $\gamma_i$ be the upper limit on the number that can be sold. In terms of options market size is often presented as (bid price)${\times}$(ask price), so here $\gamma_i{\times}\delta_i$. Then we add the constraints:
\begin{equation}
x_i \leq \delta_i~~~~\forall i \in C
\label{size1}
\end{equation}
\begin{equation}
y_i \leq \gamma_i~~~~\forall i \in C
\label{size2}
\end{equation}

\subsection{Underlying values}

To decide the values of the underlying to include in
$\Delta =[U_j~|~j=1,\ldots,N]$   first note
from Table~\ref{table1} that,  graphically,   the payoff from any option as the underlying changes is a piecewise linear curve, with a horizontal portion
  associated the fixed element ($-c_{ik}$ or $r_{ik}$, $k{=}1,2$) together with
   a linear portion with gradient  $\pm 1$ with end point at $U_j=E_i$.  In other words this payoff function is 
  continuous and piecewise linear with a break point at the exercise (strike) price.
  
 Hence the payoff function for any option portfolio composed of 
 a linear combination of calls and puts will also be continuous and piecewise linear with
respect to the value of the underlying, with any 
possible break points being associated with the exercise (strike) prices of the options chosen. For this reason we need only consider in $\Delta$ values for the underlying 
associated with the (distinct) exercise prices in $\cup_{i \in C}E_i$.  Hence $\Delta=\cup_{i \in C}E_i$.

We need to guard against negative payoffs resulting from the underlying  at expiry being outside the range $[U_1,U_N]$ considered
in $\Delta$. So here we need to consider the underlying at expiry
falling below $U_1$, or rising above $U_N$.

Consider the case of $U_1$. For values of the underlying below $U_1$
we must have an unbroken linear portfolio  payoff function (since there are no exercise prices below $U_1$). This is
 a consequence of  the continuous and piecewise linear nature of the portfolio
payoff function.

If we have a solution with $\mathcal{V}>0$ then  the payoff $\alpha_1$ at $U_1$ must also be $>0$ (from Equation~(\ref{eqvj})).
Hence to ensure that
we never have negative payoffs if the underlying falls below the lowest exercise price $U_1$ we need only adopt the constraint that the payoff at 
a single underlying value less than $U_1$ is $\geq$ the payoff at $U_1$, e.g.~the payoff at $U_1{-}1$ is $\geq$ the payoff at $U_1$, so 
the payoff below $U_1$ is non-decreasing.
This constraint prevents the portfolio payoff function below $U_1$ having a negative gradient which would (at some underlying value) lead to a loss.

 Similarly the constraint that the payoff at $U_N{+}1$ is $\geq$ the payoff at $U_N$ guards against the payoff being negative if the underlying rises above $U_N$.

Formalising this we have: $\Delta \leftarrow \Delta \cup (min[E_i~|~i \in C]{-}1) \cup (max[E_i~|~i \in C]{+}1)$, so $N \leftarrow N{+}2$ with $\alpha_1$ now begin associated with the minimum exercise price minus one and $\alpha_N$ now being associated with the maximum exercise price plus one. Then we add to the formulation
given above  the constraints:
\begin{equation}
\alpha_1 \geq \alpha_2
\label{end1}
\end{equation}
\begin{equation}
\alpha_N \geq \alpha_{N-1}
\label{end2}
\end{equation}
Essentially here, with $\Delta$ defined as above, \textbf{\emph{we have a finite set of values explicitly considered for the underlying which enable us to implicitly consider all the
(infinite) possible underlying values. }}
Obviously if a user is convinced that the entire range of underlying values at expiry is captured in the interval $[U_1,U_N]$ then there is no need to include Equations~(\ref{end1}),(\ref{end2}).

Note here that our approach as given above for considering all possible underlying values
is significantly less restrictive than that given in~\cite{beare25}. They impose the conditions that the payoff at or below the lowest exercise price
is zero and  that the payoff at or above the highest exercise price  is also zero. In terms of our notation and our $\Delta$ their approach 
reduces to setting $\alpha_1=\alpha_2=\alpha_{N-1}=\alpha_N=0$.

\subsection{Varying the budget for long and short positions}

In the above formulation the user must specify values for $\Gamma_B$, the initial cash budget allocated to buying options (so the long position) 
and $\Gamma_S$, the limit on the revenue from selling options (so the short position). 
In fact it is possible to amend the formulation to automatically decide the best possible allocation for 
$\Gamma_B$ and $\Gamma_S$ within a fixed joint budget for $(\Gamma_B  +\Gamma_S)$.

To do this and still maintain a linear formulation (having regard to computation time)  we need to drop the $z_i$ variables from the formulation. These variables were introduced to prevent  simultaneously buying and selling the same option (i.e.~to prevent $x_i>0$ and $y_i > 0$). If simultaneous buying and selling is detected in any
computational solution then we could, if desired, adjust it appropriately.

Let $\Gamma$ be the user defined limit on the total joint budget,
where $\Gamma_B$ and $\Gamma_S$ are now no longer fixed constants, but instead are variables.
 Then we amend the formulation by deleting Equation~(\ref{need3}) and adding :
\begin{equation}
\Gamma_B, \Gamma_S \geq 0
\label{add1}
\end{equation}
\begin{equation}
\Gamma_B +  \Gamma_S \leq \Gamma
\label{add2}
\end{equation}
\noindent Here Equation~(\ref{add1}) declares  $\Gamma_B$ and $\Gamma_S$ as non-negative variables and Equation~(\ref{add2}) limits their sum to be at most $\Gamma$. Equations~(\ref{eq4}),(\ref{eq3}), the equations limiting the  number of option units bought/sold  in relation to the total budget now become:
\begin{equation}
 x_i   \leq ( \Gamma_B + \Gamma_S)/c_{ik} ~~~~\forall i \in C~k{=}\tau(i)
\label{add3}
\end{equation}

\begin{equation}
 y_i   \leq  \Gamma_S/r_{ik}~~~~\forall i \in  C~k{=}\tau(i)
\label{add4}
\end{equation}

\noindent Note here that it is trivial to extend this formulation
to enable the 
the user to add, if they wish, explicit lower and upper limits on both $\Gamma_B$ and $\Gamma_S$. However if we introduce one or more lower limits then we can no longer guarantee  the existence of a solution with $\mathcal{V}=0$ and so Equations~(\ref{need1}),(\ref{need2}) are no longer valid and we must allow $\mathcal{V}$ and $[\alpha_j]$ to take both positive and negative values.

\subsection{Discussion}
The  formulation given above which enables us
to find (if possible) an option portfolio with a positive payoff is a mixed-integer linear program.
Sophisticated software packages exist for the numeric solution of such programs.
 If the solution to the formulation has value zero then it implies there does not exist any set of options (from the universe $C$ of options considered) that can lead to a profit (positive payoff).

By contrast if the solution is non-zero then it implies that the options chosen will lead to generating a \textbf{\emph{certain}} profit if they are held to expiry. 

Note here that a non-zero solution represents the \textbf{\emph{minimum certain profit}}  if the options chosen are held to expiry. The actual profit at expiry, depending upon the level of the underlying at expiry, 
might exceed this minimum value. \textbf{\emph{Note the use of the word certain here - since the above formulation deals with all possible values of the underlying we must, mathematically, be guaranteed at least this minimum profit on expiry. }}

Obviously in the real-world there are a number of caveats here. For example, whether given a non-zero (positive) solution the required option trades can be made at the prices assumed when solving the problem. This is obviously a matter relating to the technological structure adopted, e.g.~making use of automated dealing, so not requiring human intervention.

Here we regard the risk
of a counterparty defaulting on any option obligation as negligible, which seems realistic given the size of the options market and the timescale to expiry for options. For example  for S\&P~500 index options, which we consider computationally below, it is estimated that approximately \$3.5 trillion of options are traded per day, with the most commonly traded options being associated with an expiry of one month or less. 

With regard to the use of the formulation given above an obvious use is for the sole purpose of making profit, so trading simply for gain. We
can easily conceive of an automated computerised application which continually takes in current option data and solves the formulation given above. 
If no positive payoff can be achieved, then that is fine - no investment in options has been made.  After a while prices might well have changed and the
process repeats. If a positive payoff is detected then trading takes place to take advantage of it.


A less obvious use of the formulation given above  is for an institution trading in options to check whether their 
current option offerings (so calls and puts and their institution decided prices) are inadvertently exposing the institution 
to the possibility of a third party taking advantage of their offerings to generate profit for themselves.

Finally we would note here that although in this paper we only present computational results for 
 S\&P~500 index options 
 the underlying formulation applies to any asset for which  there
exist tradable options.

\section{Computational results}
\label{sec4}
\sloppy We used S\&P~500 index options
as available from
 \url{https://www.optionsdx.com/product/spx-option-chain/}. The computational results shown below are for the month of January 2023 as taken from that website and are based on end of day prices. 
We considered four cases with:
\begin{compactitem}
\item No shorting $(\Gamma_S=0)$/shorting  $(\Gamma_S>0)$
\item No size constraints/size constraints
\end{compactitem}
 For the cases with size constraints   we used the  size values as given in the data taken from
the above website to limit, for each option,  the number that could be purchased or sold at the prices given, so making use of
Equations~(\ref{size1}),(\ref{size2}).

Although there were only 19 trading days in the month we examined recall that in our approach all options
 must have the same expiry so on any particular day there may be more than one problem to solve (involving options traded that day, but with differing expiry times). 
For the month we examined there were in total 925 problems to solve  using our formulation, so an average of approximately 49 problems (49 different 
expiry times) per trading day.

Over these 925 problems the average option time to expiry was 8.59 days, with the minimum time to expiry being 0.04 days and the maximum time to expiry being 75.33 days. The  average value for $|C|$, the number of options  considered for portfolio inclusion, was 304 and the average value for
 $|\Delta|$, the number of exercises prices considered, was 154. These exercise prices spanned a wide range either side of the current underlying value.  
 For the month we examined the lowest exercise price was  on average only 16\% the current underlying value, 
 the highest exercise price was on average 157\% of the current underlying value.

We used a Windows pc with 8GB of memory, a 512GB SSD disk and an Intel Core i5-1135G7 2.4GHz processor, a multi-core pc with four cores. Cplex~\cite{cplex1210} was used as the mixed-integer solver.

Computational results can be seen in Table~\ref{table2}. We considered values for ($\Gamma_B + \Gamma_S$) of $10^3,10^4,10^5$ and $10^6$.  In terms of the underlying financial commitment each S\&P~500 index point is valued for option trading purposes as US\$100. This factor of 100 is known as the option/contract multiplier. So these values 
for ($\Gamma_B + \Gamma_S$) indicate a budget for investment of \$100K, \$1m, \$10m and \$100m. \cite{post21} refer to the first of these values as being associated with small retail investors, the other three values being associated with institutional investors.

 In Table~\ref{table2} we show the average, over the 925 problems solved, for:
\begin{compactitem}
\item the  percentage return on the money invested, so $100(\mbox{total payoff}~\mathcal{V})/ (\Gamma_B + \Gamma_S)$
\item the number of options chosen (so with $x_i > 0$ or $y_i > 0$)
\item  the percentage of problems where  there was no trade ($\mathcal{V}=0$), 
so no solution showing a positive payoff could be found
\item  the computation time in seconds.
\end{compactitem}
Note here that we adopt percentage return in this table in order to make a comparison between 
payoffs achieved using different values for the financial commitment ($\Gamma_B + \Gamma_S$).

The average for the number of days to expiry in the final column of Table~\ref{table2} is an average over the problems where there was a trade
 ($\mathcal{V}>0$). Note here that, as discussed above, the percentage return is the minimum that will certainly  be achieved
 for a given financial commitment.
  If the value of the 
underlying at expiry is more favourable for the options chosen the return achieved can exceed this minimum value.

\begin{table}[!ht]
\centering
{
\renewcommand{\tabcolsep}{1mm} 
\renewcommand{\arraystretch}{1.4} 
\begin{tabular}{|c|c|cc|c|c|c|c|c|}
\hline
&  & $\Gamma_B$ & $\Gamma_ S$ & Average &   Average &  \% no & Average & Average  \\
 &  & & &       \% return        & number of  &          trading                   & time & days to \\
  &  & & &                & options chosen &                            & (secs)  & expiry \\
 
\hline
No shorting & No size   & $10^3$   & 0  & 	0.00	&	0.04	&	98.27	& 	0.04	&	0.21	\\

& constraints & $10^4$ & 0& 	2.06	&	2.33	&	18.49	& 	0.17	&	10.43	\\

& & $10^5$ & 0 & 	2.18	&	2.49	&	18.49	& 	1.67	&	10.43	\\

& & $10^6$ & 0 & 	2.19	&	2.48	&	18.49	& 	1.20	&	10.43	\\

\cline{2-9}
   & Size   & $10^3$   & 0  &  	0.00	&	0.04	&	98.27	& 	0.04	&	0.21	\\

  & constraints  & $10^4$ & 0 & 	2.05	&	2.43	&	18.49	& 	0.18	&	10.43	\\

& & $10^5$ & 0 & 	2.14	&	3.69	&	18.49	& 	0.99	&	10.43	\\

& & $10^6$ & 0 & 	1.30	&	12.52	&	18.49	& 	0.15	&	10.43	\\

\hline
Shorting & No size   & $9{\times}10^2$   & $10^2$  &  	0.18	&	2.33	&	48.97	& 	0.54	&	1.94	\\

& constraints & $9{\times}10^3$ &  $10^3$  & 	2.05	&	6.88	&	3.57	& 	1.39	&	8.89	\\

& & $9{\times}10^4$ & $10^4$ & 	2.08	&	7.49	&	3.57	& 	3.48	&	8.89	\\

& & $9{\times}10^5$ & $10^5$ & 	2.08	&	6.94	&	3.57	& 	0.98	&	8.89	\\

\cline{2-9}
 & Size   & $9{\times}10^2$   & $10^2$  & 	0.18	&	2.36	&	49.08	& 	0.61	&	1.94	\\

& constraints & $9{\times}10^3$ &  $10^3$  & 	2.04	&	8.13	&	3.68	& 	1.49	&	8.90	\\

& & $9{\times}10^4$ & $10^4$ & 	2.02	&	18.49	&	3.68	& 	11.77	&	8.90	\\

& & $9{\times}10^5$ & $10^5$ 	& 	1.50	&	53.77	&	3.68	& 	11.14	&	8.90	\\

\cline{2-9}
 & No size   & 0 & $10^3$  & 	0.27	&	2.14	&	62.49	& 	0.44	&	1.25	\\

& constraints & 0 & $10^4$   & 	0.28	&	2.77	&	60.65	& 	0.55	&	1.39	\\

& & 0 &  $10^5$	& 	0.28	&	3.34	&	60.43	& 	1.03	&	1.50	\\

& & 0 & $10^6$ 	& 	0.28	&	3.23	&	60.43	& 	0.65	&	1.50	\\

\cline{2-9}
 & Size   & 0 & $10^3$   & 	0.24	&	2.86	&	62.59	& 	0.50	&	1.25	\\

& constraints & 0 & $10^4$   & 	0.12	&	8.18	&	61.30	& 	4.52	&	1.37	\\

& & 0 &  $10^5$	& 	0.02	&	10.96	&	61.30	& 	8.99	&	1.37	\\

& & 0 & $10^6$ & 	0.00	&	11.93	&	61.30	& 	10.76	&	1.37	\\

\hline
\end{tabular}
}
\caption{Computational results, $\Gamma_B$ and $\Gamma_S$ fixed}
\label{table2}
\end{table}

In Table~\ref{table2} for the no shorting case (so $\Gamma_S=0$) we can see that as we would expect, imposing size constraints results in a  
decrease in return, an increase in the number of options chosen and an increase in the time taken. 

It is clear from that table that, for no shorting
with
size constraints, moving from a budget of $\Gamma_B=10^5$ to $\Gamma_B=10^6$ reduces the average return. 
The reason for this is that providing a larger 
budget enables more trades to be made, as signified by the significant increase in the number of options chosen, and too increases the payoff.
However since for comparison purposes in Table~\ref{table2} we divide the payoff by the investment ($\Gamma_B + \Gamma_S$) a larger budget 
 can result in the percentage
return decreasing, despite the increase in payoff.

In Table~\ref{table2} when we have  shorting  
we  consider two situations, the first with the budget for shorting being 10\% of the overall budget, but with 
the overall budgets being the same as for  no shorting. 

 In this first situation note that  with a total
budget of $\Gamma_B + \Gamma_S = 10^6$  we see a significant increase in the 
number of options chosen and a significant decrease in return, so the same behaviour as observed for no shorting with size constraints.

The second shorting situation considered in Table~\ref{table2} has $\Gamma_B=0$, so we have no initial budget for purchases, rather any purchases must be financed by short selling,  so making use of the $\Gamma_S$ budget. The returns seen for this situation are markedly reduced from the first shorting situation considered, with around 60\% of days having no trade. 

It is clear from Table~\ref{table2} that, at least for the data we considered, having no initial budget for purchases and instead relying on generating funds from short selling is inferior to the other cases considered for  $\Gamma_B + \Gamma_S \geq 10^4$.

Table~\ref{table3} shows the results when we allow $\Gamma_B$ and  $\Gamma_S$ to be variable
and has the same format as Table~\ref{table2}.
Here we examine the same values for $\Gamma = \Gamma_B + \Gamma_S$ as we used in Table~\ref{table2}.


As we would expect, since here we now have more flexibility in terms of $\Gamma_B$ and $\Gamma_S$,
all of the returns given in Table~\ref{table3}
dominate the corresponding entries in Table~\ref{table2}.

\begin{table}[!ht]
\centering
{
\renewcommand{\tabcolsep}{1mm} 
\renewcommand{\arraystretch}{1.4} 
\begin{tabular}{|c|c|c|c|c|c|c|}
\hline
& $\Gamma$  & Average &   Average &  \% no & Average & Average  \\
 & &       \% return          & number of  &          trading                   & time & days to  \\
& &                 & options chosen &                            & (secs)  & expiry   \\
 
\hline
No size   & $10^3$  & 	1.01	&	3.92	&	17.84	& 	0.77	&	5.77		\\

constraints & $10^4$ & 	2.43	&	5.97	&	3.57	& 	0.71	&	8.89		\\

& $10^5$ & 	2.47	&	6.55	&	3.57	& 	4.53	&	8.89		\\

& $10^6$	& 	2.47	&	6.55	&	3.57	& 	4.53	&	8.89		\\

\hline

Size   & $10^3$   	& 	0.98	&	4.55	&	17.95	& 	0.68	&	5.78	\\

constraints & $10^4$	& 	2.28	&	10.53	&	3.57	& 	1.53	&	8.89		\\

& $10^5$ 	& 	2.22	&	17.35	&	3.57	& 	7.64	&	8.89		\\

& $10^6$	& 	1.59	&	48.98	&	3.68	& 	7.20	&	8.90		\\

\hline
\end{tabular}
}
\caption{Computational results, $\Gamma_B$ and $\Gamma_S$ variable}
\label{table3}
\end{table}

Of particular note here with regard to both Table~\ref{table2} and Table~\ref{table3} is the 
low computation time and the average percentage return. Given that we 
might expect prices to change relatively rapidly in the market achieving a solution quickly will be of 
importance. Also note that the percentage return relates to the option portfolio held to expiry. 
The returns achieved  over the relatively short timescale involved (days to expiry, as detailed in Table~\ref{table2} and Table~\ref{table3})  seems to illustrate  
the potential of the approach given in this paper – at least on the data we examined.

\section{Conclusions}
\label{sec5}

In this paper we presented  a mixed-integer linear program for the problem of finding a portfolio of options that will give a guaranteed minimum payoff irrespective of the value of the underlying at expiry. In our approach all of the options considered must have the same expiry.

We discussed the use that might be made of our approach and gave computational results related to constructing option portfolios made up from S\&P~500 index options. These results indicated that, at least for the data we examined, there is potential for achieving positive payoffs.

\FloatBarrier
\linespread{1}
\small \normalsize 
\bibliographystyle{apalike}
\bibliography{paper}

\end{document}